\documentclass[final,12pt,nonatbib]{elsarticle}
{
\makeatletter
\let\c@author\relax
\makeatother
}

\usepackage{xcolor,soul,cancel}
\usepackage[color]{showkeys}

\usepackage{mathtools}
\definecolor{refkey}{rgb}{0,1,1}
\definecolor{labelkey}{rgb}{1,0,0}
\usepackage{amssymb}
\usepackage{amsmath}
\usepackage{amsthm}
\usepackage{graphicx}
\usepackage{float}
\usepackage{xcolor}
\usepackage{cases}
\usepackage{color}
\usepackage{hyperref}
\hypersetup{
    colorlinks=true, 
    linkcolor=blue, 
    urlcolor=red, 
    linktoc=all 
}
\journal{arXiv}

\newtheorem{thm}{Theorem}
\newtheorem{lem}{Lemma}
\newtheorem{cor}{Corollary}
\newtheorem{prop}[thm]{Proposition}

\numberwithin{equation}{section}
\newcommand{\floor}[1]{\lfloor #1 \rfloor}
\newcommand{\eq} [1] {\begin{equation}\label{#1}\quad}
\newcommand{\en} {\end{equation}}

\newcommand{\norm}[1]{\left\Vert#1\right\Vert}
\newcommand{\abs}[1]{\left\vert#1\right\vert}

\newcommand{\C}{\mathbb C}
\newcommand{\M}{{\bf M}}
\newcommand{\R}{\mathbb R}

\newcommand{\im}{\operatorname{Im}}

\newcommand{\conv}{\operatorname{conv}}
\newcommand{\re}{\operatorname{Re}}

\newcommand{\Rec}{\operatorname{\bf RM}}

\begin{document}

\begin{frontmatter}

\title{Low-dimensional reciprocal matrices with elliptical components of their Kippenhahn curves\tnoteref{support}}


\author[nyuad]{Muyan Jiang}
\ead{mj2259@nyu.edu, jimmymuyan@163.com}
\author[nyuad]{Ilya M. Spitkovsky}
\address[nyuad]{Division of Science and Mathematics, New York  University Abu Dhabi (NYUAD), Saadiyat Island,
P.O. Box 129188 Abu Dhabi, United Arab Emirates}
\ead{ims2@nyu.edu, ilya@math.wm.edu, imspitkovsky@gmail.com}
\tnotetext[support]{This work is based on the capstone project of the first author [MJ] supervised by the second author [IMS]. The latter was supported in part by Faculty Research funding from the Division of Science and Mathematics, New York University Abu Dhabi.}

\begin{abstract}
By definition, reciprocal matrices are tridiagonal $n$-by-$n$ matrices $A$ with constant main diagonal and such that $a_{i,i+1}a_{i+1,i}=1$ for $i=1,\ldots,n-1$. For $n\leq 6$, we establish criteria under which the numerical range generating curves (also called Kippenhahn curves) of such matrices consist of elliptical components only. As a corollary, we also provide a complete description of higher rank numerical ranges when the criteria are met. 
\end{abstract}


\end{frontmatter}

\section{Introduction} 
Denote by $\M_n$ the algebra of all $n$-by-$n$ matrices with complex entries. For any $A\in\M_n$, we will use the standard notation 
$\re A=(A+A^*)/2,\ \im A= (A-A^*)/(2i)$. Letting \eq{eig} \lambda_1(\theta)\geq\ldots\geq\lambda_n(\theta) \en stand for the eigenvalues (counting their multiplicities) of the hermitian matrix $\re(e^{i\theta}A)$, introduce the family of lines  \eq{l} \ell_j(\theta) =e^{-i\theta}(\lambda_j(\theta)+i\R), \quad \theta\in[0,2\pi),\ j=1,\ldots,n,\en
and the respective family of half-planes $H_j(\theta)=\{ z\in\C\colon \re(e^{i\theta}z)\leq \lambda_j(\theta)\}$ bounded by $\ell_j(\theta)$. For every $k=1,\ldots,n$ the intersection \eq{rkf} \Lambda_k(A)=\cap_{\theta\in[0,2\pi)} H_k(\theta)\en  is nothing but the {\em rank-$k$ numerical range} of $A$ defined as 
\eq{rk} \Lambda_k(A)=\{\lambda\in\C\colon PAP=\lambda P\text{ for some rank-}k\text{ orthogonal projection }P\}.\en 
Obviously, 
	\eq{chain}\Lambda_1(A)\supseteq \Lambda_2(A) \supseteq \ldots \supseteq \Lambda_n(A). \en 
Note that $\Lambda_1(A)$ is the regular {\em numerical range} $W(A)$ of $A$, and the respective formula \eqref{rkf} is an immediate consequence of its convexity (the celebrated Toeplitz-Hausdorff theorem, see e.g. \cite{GusRa}). For $k>1$, however, \eqref{rkf} is a more recent result \cite{LiSze08}. In particular, all $\Lambda_k(A)$ are convex \cite{LiSze08,Woe08}. 

Due to \eqref{rkf}, all the sets $\Lambda_k(A)$ are defined completely by the envelope $C(A)$ of the family \eqref{l}. In particular, $W(A)$ is the convex hull of $C(A)$, as was already observed in \cite{Ki} (see \cite{Ki08} for the English translation). Called the ``boundary generating curve'' in \cite{Ki}, in the recent years $C(A)$ was renamed as {\em Kippenhahn curve} --- the term which we will be using in what follows. 

Along the same lines (no pun intended), the family \eqref{l} is defined by \eqref{eig}, and thus by the characteristic polynomial 
\eq{pol} P(\lambda,\theta)=\det(\re(e^{i\theta}A)-\lambda I)\en of $\re(e^{i\theta}A)$; the latter is therefore called the {\em Kippenhahn polynomial} of $A$.

A vast number of papers is devoted to singling out classes of matrices $A$ with $W(A)$ being an elliptical disk, the convex hull of several ellipses, or at least having an elliptical arc as part of its boundary $\partial W(A)$. Similar questions can be asked about higher rank numerical ranges. In terms of Kippenhahn curves, this boils down to the following: when does $C(A)$ contain an ellipse $E$ as one of its components (or, more restrictively, consists of several ellipses)?  

In \cite{BPSV} this question was addressed for so called {\em reciprocal matrices} of small (up to $n=6$) sizes, provided that the ellipse $E$ is centered at the origin. This condition on $E$ is lifted here. Also, the higher rank numerical ranges of the respective matrices $A$ are described.

\section{General observations} 

The classical Elliptical range theorem claims that the numerical range $W(A)$ of $A\in\M_2$ is an elliptical disk with the foci at the eigenvalues $\zeta_1,\zeta_2$ of $A$ and the minor axis of the length $2c:=\sqrt{\norm{A}_{\rm Fr}^2-\abs{\zeta_1}^2-\abs{\zeta_2}^2}$ 
(here $\norm{A}_{\rm Fr}$ denotes the Frobenius norm of $A$). Note that $c=0$ if and only if $A$ is normal, and then $W(A)$ degenerates into the line segment $[\zeta_1,\zeta_2]$. 

The $n=2$ case is exceptional in a sense that $W(A)$ defines $A$ up to a unitary similarity and, when put in an upper-triangular form, it becomes $A=\begin{bmatrix}\zeta_1 & 2c \\ 0 & \zeta_2\end{bmatrix}$. 
A direct computation shows that the Kippenhahn polynomial of this matrix is \eq{pole}P_{E} = (\lambda -p \cos (\theta )+q \sin (\theta ))^2-(x \cos (2 \theta )-y \sin (2 \theta )+z), \en 
where 
\eq{xyz}
\begin{aligned}
	    p & = \frac{1}{2}\re(\zeta_1+\zeta_2), \\
	    q & = \frac{1}{2}\im(\zeta_1+\zeta_2), \\ 
	    x & = \frac{1}{8} \re\left( (\zeta_1-\zeta_2)^2\right), \\
		y & = \frac{1}{4}  \re (\zeta_1-\zeta_2) \im (\zeta_1-\zeta_2),\\
		z & = \frac{1}{8}  \abs{\zeta_1-\zeta_2}^2+c^2.
\end{aligned}\en 
The third and the fourth equations of \eqref{xyz} imply that $\sqrt{x^2+y^2}=\frac{1}{8}\abs{\zeta_1-\zeta_2}$ and so, invoking the fifth: \eq{c}c^2=z-\sqrt{x^2+y^2}. \en 
In particular, \eq{zin} z\geq \sqrt{x^2+y^2},\en  and the equality holds if and only if the quadratic polynomial \eqref{pole} factors into the product of two linear ones. The ellipse $C(A)$ then collapses into the doubleton $\{\zeta_1,\zeta_2\}$ of its foci.  

Turning to matrices of arbitrary size, we therefore have 
\begin{prop}\label{th:qfac}Let $A\in\M_n$. Then $C(A)$ contains an ellipse $E$ if and only if the Kippenhahn polynomial \eqref{pol} of $A$ is divisible by \eqref{pole} for some $p,q,x,y\in\R$ and $z\in\R_+$ satisfying \eqref{zin}. If this condition holds, then $E$ is centered at $(p,q)$, the length of its minor half-axis $c$ is given by \eqref{c}, and the difference between the foci is $8\sqrt{x^2+y^2}\exp(i(\tan^{-1} (y/x)/2))$. \end{prop}
Proposition~\ref{th:qfac} is not new; its particular case corresponding to the ellipticity criterion for $W(A)$ in terms of an explicit formula for $\lambda_1(\theta)$ is \cite[Theorem 1]{ChienHung}. In slightly different notation, and without specifying the minor axis and foci of $E$, it was also used in \cite{GeS2}. We chose to state it in a detailed form for convenience of reference. 

Following \cite{BPSV}, we say that a matrix $A=(a_{ij})_{i,j=1}^n$ is {\em reciprocal} if it is tridiagonal with zero main diagonal, i.e. $a_{ij}=0$ whenever $\abs{i-j}\neq 1$, while \eq{rec} a_{i,i+1}a_{i+1,i}=1 \text{ for } i=1,\ldots,n-1.\en 

In what follows, we will denote by $\Rec_n$ the class of all reciprocal matrices $A\in\M_n$. 

Along with such $A$, all the matrices $\re(e^{i\theta}A)$ are tridiagonal with the zero main diagonal. This observation was used in \cite[Proposition 4]{BPSV} to conclude that Kippenhahn polynomials of reciprocal matrices can be written as \[ P_n(\zeta,\tau)=\zeta^k+\sum_{j=0}^{k-1}p_j(\tau)\zeta^j,\]
premultiplied by $-\lambda$ if $n$ is odd. Here $\zeta=\lambda^2$, $k=\floor{n/2}$, and $p_j$ are polynomias of degree $k-j$ in $\tau:=\cos(2\theta)$ with the coefficients depending only on \eq{Aj} A_j:= \frac{\abs{a_{j,j+1}}^2+\abs{a_{j+1,j}}^2}{2},\quad j=1,\ldots,n-1.\en  Respectively, $C(A)$ is symmetric about both coordinate axes \cite[Corollary 2]{BPSV}, for odd $n$ including the origin as one of its components. 

In this paper, we find it more convenient to represent $P_n$ as 
\eq{polrec} P_n(\zeta,\rho)=\zeta^k+\sum_{j=0}^{k-1}\widetilde{p_j}(\rho)\zeta^j,\en
with the variable $\tau$ replaced by $\rho:=\cos^2\theta\left(=(1+\tau)/2\right)$ and the coefficients of $\widetilde{p_j}$ expressed in terms of 
\eq{xij} \xi_j:= \frac{\abs{a_{j,j+1}}-\abs{a_{j+1,j}})^2}{4} \ \left(=\frac{A_j-1}{2}\right),\quad j=1,\ldots,n-1.\en 
Note that $\xi_j\geq 0$, and the equality holds if and only if $\abs{a_{j,j+1}}$ (equivalently, $A_j$, or $\abs{a_{j+1,j}}$) is equal to one.
\begin{prop}\label{th:real}Suppose the Kippenhahn curve $C(A)$ of a reciprocal matrix $A$ contains an ellipse $E$. Then {\em (i)} the foci of $E$ are real, and {\em (ii)} either $E$ is centered at the origin, or its reflection $-E$ is also contained in $C(A)$. \end{prop}
\begin{proof}Part (ii) follows immediately from the symmetry of $C(A)$. To prove (i), recall that for any square matrix the foci of its Kippenhahn curve coincide with its spectrum \cite[Theorem 11]{Ki08}. If $A$ is reciprocal, it is similar (not unitarily similar!) to the Toeplitz tridiagonal matrix $T$, as was observed in \cite[Proposition 5]{BPSV} by the anonymous referee's suggestion. Consequently, for any $A\in\Rec_n$ its spectrum 
\eq{sA} \sigma(A)=\sigma(T)=\left\{ 2\cos\frac{j\pi}{n+1}\colon j=1,\ldots,n)\right\} \en
is real. The foci of $E$ lie in the set \eqref{sA}, and as such are also real. 
\end{proof} 

\begin{cor}\label{th:bigk}For any $A\in\Rec_n$, rank-$k$ numerical ranges with $k>(n+1)/2$ are empty. \end{cor}
\begin{proof}Indeed, $\Lambda_k(A)$ for $k>n/2$ can only be a singleton coinciding with an eigenvalue of $A$ of multiplicity at least $2k-n$ \cite[Proposition 2.2]{ChoKriZy06}, and \eqref{sA} implies that all the eigenvalues of $A$ are simple. \end{proof} 
According to Proposition~\ref{th:real}(i), for an ellipse $E\subset C(A)$ in case of reciprocal $A$ we have $q=y=0$, with the rest of formulas \eqref{xyz} simplifying to 
\eq{xyz1} p=\frac{1}{2}(\zeta_1+\zeta_2),\ x=\frac{1}{8}(\zeta_1-\zeta_2)^2,\ z= x+c^2. \en 
Proposition~\ref{th:qfac} for reciprocal matrices can therefore be recast as follows.
\begin{prop}\label{th:facrec}Let $A\in\Rec_n$. Then $C(A)$ contains an ellipse with the half-axes of length $c$ and $\sqrt{c^2+X^2}$ if and only if polynomial \eqref{polrec} is divisible either {\rm (a)} by $\zeta-(X^2\rho+c^2)$, or {\rm (b)} by  
\eq{2el} \zeta^2-2\zeta\left((X^2+p^2)\rho+c^2\right)+\left((X^2-p^2)\rho+c^2\right)^2. \en 
\end{prop} 
Case (a) corresponds to an ellipse centered at the origin. In case (b) there are two symmetric ellipses $\pm E$ contained in $C(A)$ by Proposition~\ref{th:real}(ii); $\pm p$ denote their centers. Formula \eqref{2el} is obtained by multiplying out the respective polynomials 
$P_{\pm E}=(\lambda\mp p\cos\theta)^2-(x\cos(2\theta)+z)$; in both cases we relabeled $2x=:X^2$.   

There is another important consequence of $\re(e^{i\theta}A)$ being tridiagonal hermitian matrices. According to \cite[Corollary 7]{BS041}, such matrices can have a repeated eigenvalue only if they split into the direct sum of two blocks sharing this eigenvalue. Due to \eqref{rec}, in our setting this can only happen if $\theta=\pi/2\mod\pi$ and $A_k=1$ (equivalently: $\xi_k=0$) for some $k$.
\begin{prop}\label{th:multan}For a reciprocal matrix $A$, the multiple tangent lines of $C(A)$ can only be horizontal. Such lines materialize if and only if $\xi_k=0$ for some $k$, and the spectra of the left upper $k$-by-$k$ and the right lower $(n-k)$-by-$(n-k)$ block
of the matrix $K=\im A$ overlap. 
\end{prop} 
\begin{proof}We just need to recall that tangent lines to $C(A)$ form the family \eqref{l}.  \end{proof}
Due to the symmetry of $C(A)$, its multiple tangent lines, if any, come in pairs symetric with respect to the abscissa axis. Their ordinates are the multiple eigenvalues of $\im A$. 

Note that the presence of an ellipse  $E\subset C(A)$ centered at $p\neq 0$ implies the existence of lines tangent to both $E$ and $-E$. Proposition~\ref{th:multan} thus implies 
\begin{cor} \label{th:Aj1}Let $A\in\Rec_n$ be such that $C(A)$ contains an ellipse $E$ with its center different from the origin. Then {\em (i)} $E\cap (-E)\neq\emptyset$ and {\em (ii)} $\abs{a_{k,k+1}}=1$ (equivalently: $\xi_k=0$) for some $k\in\{2,\ldots,n-2\}$, while the upper left $k$-by-$k$ block of $\im A$ has a non-zero eigenvalue in common with its lower right $(n-k)$-by-$(n-k)$ block. \end{cor}  
\begin{proof}Being congruent, the ellipses $E$ and $-E$ either intersect or lie outside of each other. In the latter case, they would have a non-horizontal tangent in common which would contradict Proposition~\ref{th:multan}. This proves (i). 
		
Part (ii) also follows from Proposition~\ref{th:multan} as soon as we observe that a shared eigenvalue equal to zero corresponds to $E$ degenerating into the pair of its foci. If $k=1$ or $n-1$, one of the blocks of $K$ is a one-dimensional $\{0\}$ and, as such, cannot generate a non-zero eigenvalue. \end{proof} 

In contrast with Corollary~\ref{th:Aj1}(i), the concentric elliptical components of $C(A)$, if any, have to be nested. 
\begin{cor}\label{th:concen}Suppose the Kippenhahn curve of a reciprocal matrix $A$ contains two concentric ellipses. Then one of them has to lie inside the other. \end{cor}
\begin{proof}The reason is exactly the same as in the proof of the previous result: the absence of non-horizontal multiple tangents of $C(A)$. \end{proof} 
In conclusion, note that reversing the order of rows and columns of a matrix is a unitary similarity which preserves the reciprocal structure, while switching the super- and subdiagonal and reversing the order of their entries. Therefore, the Kippenhahn polynomials and curves of reciprocal matrices are invariant under the transformation \eq{flip} \xi_j\longleftrightarrow \xi_{n-j+1}, \quad j=1,\ldots, n-1.\en 
This simple observation will prove itself useful in the next sections.
\section{Reciprocal 4-by-4 and 5-by-5 matrices} \label{s:45}
Let us start with $A\in\M_4$. The Kippenhahn polynomial \eqref{pol} of $A$ then has degree four in $\lambda$. So, if a quadratic polynomial  can be factored out, the remaining multiple is also quadratic.  In the language of Kippenhahn curves it means that if $C(A)$ contains an elliptical component, it actually consists of two ellipses.  

If $A\in\Rec_4$, Proposition~\ref{th:real}(ii) implies in addition that the two ellipses in question are either both centered at the origin, or are reflections of each other through the origin (equivalently: across the ordinate axis).   

According to \eqref{sA} with $n=4$: 
\eq{sA4} \sigma(A)=\{\phi,\phi^{-1},-\phi,-\phi^{-1}\},\en 
where $\phi=(\sqrt{5}+1)/2$ is the golden ratio, while \eqref{polrec} takes the form 
\eq{char4} P_4(\zeta,\rho)=\zeta^2-(\xi_1+\xi_2+\xi_3+3\rho)\zeta+(\xi_1+\rho)(\xi_3+\rho). \en 
In particular, the characteristic polynomial of $\im A$ is 
\eq{charK4} P_4(\lambda^2,0)=\lambda^4-\lambda^2(\xi_1+\xi_2+\xi_3)+\xi_1\xi_3.  \en
\begin{thm}\label{th:el4}The Kippenhahn curve $C(A)$ of $A\in\Rec_4$ consists of two ellipses if and only if either
\eq{con4} \xi_2=\phi\xi_1-\phi^{-1}\xi_3 \text{ or } \xi_2=\phi\xi_3-\phi^{-1}\xi_1,\en
with at least one of $\xi_j$ being different from zero, or 
\eq{noncon4} \xi_2=0, \quad \xi_1=\xi_3\neq 0. \en 
\end{thm} 
Case \eqref{con4} corresponds to ellipses centered at the origin, and the criterion (up to the notational change from $A_j$ to $\xi_j$) is \cite[Theorem 8]{BPSV}. According to \eqref{sA4}, one of the ellipses (say, $E_1$) has foci $\pm\phi$ while the foci of $E_2$ then are $\pm\phi^{-1}$. The lengths of the minor half-axes are determined by the positive roots of \eqref{charK4}, and direct computations yield $\sqrt{\xi_1}\phi,\ \sqrt{\xi_3}\phi^{-1}$ or $\sqrt{\xi_3}\phi,\ \sqrt{\xi_1}\phi^{-1}$, depending on which of the equalities holds in \eqref{con4}. 

Note that the situation \eq{ext} \xi_1=0,\ \xi_2=\phi\xi_3\neq 0 \text{ or } \xi_3=0,\ \xi_2=\phi\xi_1\neq 0\en  falls under \eqref{con4} and is formally treated as two concentric ellipses, in spite of the fact that the inner one degenerates into the pair of its foci. Also, under condition \eqref{con4} at least one of $\xi_j$ being different from zero implies that only one of them can actually equal zero. 

Criterion \eqref{noncon4} (once again, up to the notational change from $A_j$ to $\xi_j$) was derived in \cite[Theorem 6.1]{GeS2} as a corollary of a more general Theorem~5.1 on 4-by-4 matrices with scalar diagonal blocks. Here is a streamlined reasoning, specific for the case at hand. 

{\sl Necessity.} Suppose $C(A)=E_1\cup E_2$ with $E_{1,2}$ being symmetric images of each other. According to Corollary~\ref{th:Aj1}, $\xi_2=0$ and $\sigma(\im A)$ consists of two opposite eigenvalues, each of multiplicity 2. From \eqref{charK4} we conclude that $\xi_1=\xi_3$.

{\sl Sufficiency.} Let \eqref{noncon4} hold. Denoting the common value of $\xi_{1,3}$ by $c^2$, it is easy to see that \eqref{char4} can be rewritten as \eqref{2el} with $p=1/2$ and $X=\sqrt{5}/2$. Proposition~\ref{th:facrec} then implies that $E_1$ and $E_2$ are congruent ellipses with the foci $\phi,-\phi^{-1}$ and  $-\phi,\phi^{-1}$, respectively, and $c$ as the length of their minor half-axes. 
\hfill $\square$

Theorem~\ref{th:el4} immediately implies that $W(A)$ is the elliptical disk bounded by $E_1$ if \eqref{con4} holds, and $\conv\{E_1\cup E_2\}$ under condition \eqref{noncon4}. These results are also from \cite{BPSV} and \cite{GeS2}, respectively. 

Formulas \eqref{rkf} allow us to move further along the chain \eqref{chain}.
\begin{cor}\label{th:hk} Let $A\in\Rec_4$ satisfy \eqref{con4} or \eqref{noncon4}. Then $\Lambda_2(A)$ is the elliptical disk bounded by the inner ellipse $E_2$ in case \eqref{con4}, and the intersection of the elliptical disks bounded by $E_1$ and $E_2$ in case \eqref{noncon4}. \end{cor}  
Note that $\Lambda_2(A)$ is the line segment $[-\phi^{-1},\phi^{-1}]$ in the extreme subcase \eqref{ext} of \eqref{con4}. Also, $\Lambda_2(A)\neq\emptyset$. This is obvious under \eqref{con4}, and follows from Corollary~\ref{th:Aj1}(ii) under \eqref{noncon4}. 

On the other hand, for $A\in\Rec_4$ we have $\Lambda_3(A)=\Lambda_4(A)=\emptyset$ by Corollary~\ref{th:bigk}.

We now move to $n=5$. The Kippenhahn polynomial of $A\in\Rec_5$ is a product of $-\lambda$ by 
\eq{char5}\begin{aligned} P_5(\zeta,\rho)   = \zeta^2 & -(\xi_1+\xi_2+\xi_3+\xi_4+4\rho)\zeta +(\xi_1+\rho)(\xi_3+\rho)
\\ &  +(\xi_1+\rho)(\xi_4+\rho)+(\xi_2+\rho)(\xi_4+\rho).\end{aligned} \en
Qualitatevely, \eqref{char5} is similar to \eqref{char4}, implying that in our current setting again $C(A)$ contains two ellipses as soon as it contains at least one. Moreover, these ellipses are either both centered at the origin (which for $n=5$ is the third component of $C(A)$), or are reflections of each other through the ordinate axis. The foci of these ellipses are located at the non-zero points of $\sigma(A)$, which are $\pm\sqrt{3},\ \pm 1$. 

According to \cite[Theorem 9]{BPSV}, for two concentric ellipses to materialize it is necessary and sufficient that 
\eq{con5}\xi_1=\xi_4\text{ or } \xi_1-\xi_4=2(\xi_3-\xi_2), \en
with not all of $\xi_j$ equal zero. Note that due to \eqref{con5} at least two of $\xi_j$ are then different from zero. Computing the roots of \eqref{char5} with $\rho=0$, we find that under either of equalities \eqref{con5} the lengths of the minor half-axes of the ellipses are 
\[ \sqrt{(\xi_1+\xi_4)/2} \text{ and  }  \sqrt{\xi_2+\xi_3+(\xi_1+\xi_4)/2}. \]
By Corollary~\ref{th:concen}, they correspond to the ellipses with their foci at $\pm 1$ and $\pm\sqrt{3}$, respectively. The inner ellipse degenerates if and only if $\xi_1=\xi_4=0$. 

Passing to the case of displaced ellipses, let us first derive the criterion for $C(A)$ to have multiple tangent lines. 
\begin{lem} \label{th:hor5}A matrix $A\in\Rec_5$ has the Kippenhahn curve with multiple tangeint lines if and only if 
\eq{mt5} \xi_1+\xi_2=\xi_3+\xi_4\neq 0\text{ while } \xi_2\xi_3=0. \en 
 \end{lem} 
\begin{proof}By Proposition~\ref{th:multan}, only horizontal multiple tangent line are possible, and $\xi_2\xi_3=0$ is a  necessary condition for them to materialize. Due to the invariance of \eqref{mt5} under \eqref{flip}, we may without loss of generality suppose that $\xi_2=0$ and concentrate on showing that then \eq{xi134} \xi_1=\xi_3+\xi_4\ (\neq 0) \en  is the desired criterion. 

Condition $\xi_2=0$ implies that $\im A =B_1\oplus B_2$, where 
\eq{b} B_1=\begin{bmatrix}0 & \eta_1 \\ \overline{\eta_1} & 0\end{bmatrix}, \quad B_2=\begin{bmatrix}0 & \eta_3 & 0 \\ 
\overline{\eta_3} & 0 &  \eta_4 \\ 0  & \overline{\eta_4} & 0\end{bmatrix},\text{ and } \eta_j=\frac{a_j-1/\overline{a_j}}{2i}.\en
The remaining part of the requirements of Proposition~\ref{th:multan} is that the spectra of $B_1$ and $B_2$ overlap; in other words, 
$\abs{\eta_1}^2=\abs{\eta_3}^2+\abs{\eta_4}^2$. Since $\abs{\eta_j}^2=\xi_j$, the necessity is proven. The case $\xi_1=0$ has to be excluded, because otherwise $C(A)$ degenerates into $\sigma(A)$. 

The sufficiency can be demonstrated via a straightforward computation. Since the blocks $B_1,B_2$ are adjacent, it also follows from  
the result of \cite[Theorem 10]{BS041} applicable to arbitrary tridiagonal matrices. 
\end{proof} 
Note that the multiple eigenvalues of $\im A$ are in fact its maximal and minimal ones. So, the horizontal tangent lines of $C(A)$ are the supporting lines of $W(A)$. Lemma~\ref{th:hor5} therefore delivers the criterion for the numerical range of a reciprocal 5-by-5 matrix to have flat portions on its boundary.
\begin{thm}\label{th:el5}The Kippenhahn curve of $A\in\Rec_5$ contains two non-concentric ellipses if and only if 
\eq{noncon5} \xi_2\xi_3=0,\ \xi_2+\xi_3\neq 0,\ \xi_1=\frac{\sqrt{3}}{2}(\xi_2+\xi_3)+\xi_3,\ 
\xi_4=\frac{\sqrt{3}}{2}(\xi_2+\xi_3)+\xi_2. \en
If this is the case, then $C(A)=E_1\cup E_2\cup\{0\}$, with $\sqrt{3},\ -1$ and $-\sqrt{3},\ 1$ being the foci of $E_1$ and $E_2$, respectively, and the length of the minor axis for both equals $(1+\sqrt{3})\sqrt{\xi_2+\xi_3}$.   \end{thm} 
\begin{proof}Conditions of Lemma~\ref{th:hor5} are necessary for $C(A)$ to contain non-concentric ellipses, so we may suppose that \eqref{mt5} holds. Invoking the symmetry under \eqref{flip}, we may suppose further that $\xi_2=0$ and \eqref{xi134} holds, while \eqref{noncon5} simplifies to
\eq{noncon51} \xi_2=0,\ \xi_3\neq 0,\ \xi_1=\frac{2+\sqrt{3}}{2}\xi_3,\ \xi_4=\frac{\sqrt{3}}{2}\xi_3. \en		
		
Polynomial \eqref{char5} in turn takes the form \[ \zeta^2-2(\xi_1+2\rho)\zeta+(\xi_1+\rho)(\xi_1+2\rho)+\rho(\xi_4+\rho),\]
which equals \eqref{2el} if and only if \eq{xic} \xi_1=c^2,\en  
\eq{xip} X^2+p^2=2,\quad (X^2-p^2)^2=3, \en and \eq{xi4} 2c^2(X^2-p^2)=3\xi_1+\xi_4. \en From \eqref{xi4} it follows that $X^2-p^2>0$, 
and so \eqref{xip} is equivalent to \[ X=\frac{\sqrt{3}+1}{2},\quad p=\frac{\sqrt{3}-1}{2}. \]  Solving \eqref{xi134}, \eqref{xi4} for $\xi_1,\xi_4$ gives \eqref{noncon51}. It remains to invoke Proposition~\ref{th:facrec} to prove the criterion. The centers of $E_{1,2}$ are located at $\pm p$, which uniquely determines their foci, and \eqref{xic} combined with \eqref{noncon51} yields the formula for the length $2c$ of the minor axis.
\end{proof} 
Description of $\Lambda_2(A)$ for $n=5$ is very similar to Corollary~\ref{th:hk}. Namely:
\begin{cor}\label{th:hk5}Let $A\in\Rec_5$ satisfy \eqref{con5} or \eqref{noncon5}. Then $\Lambda_2(A)$ is the elliptical disk with the foci $\pm 1$ and the minor axis of the length $\sqrt{2(\xi_1+\xi_4)}$ in the former case, and the intersection of the elliptical disks bounded by $E_1, E_2$ described in Theorem~\ref{th:el5} in the latter. \end{cor}

As opposed to the case $n=4$, for $n=5$ we have $\Lambda_3(A)=\{0\}$, not the empty set. On the other hand, $\Lambda_4(A)$ is still empty by Corollary~\ref{th:bigk}. 

Finally, observe that \eqref{xi134} may hold while \eqref{noncon5} does not. This means that the numerical range $W(A)$ of a reciprocal 5-by-5 matrix may have flat portions on its boundary while the smooth arcs connecting these (horizontal) flat portions are not elliptical. 
The proof of the criterion \eqref{noncon4} outlined on p.~8 shows that for $n=4$ this is an impossibility. 

To illustrate, consider $A\in\Rec_5$ with \eq{xip} \xi_1=\xi_3=0.5,\ \xi_2=\xi_4=0.\en  Conditions \eqref{mt5} hold for this $A$, while \eqref{noncon51} fail. The Kippenhahn curve $C(A)$ consists of the two non-elliptical (``drop-shaped'') symmetric curves $D$, $-D$, and the origin. Note that qualitatively the structure of the rank-$k$ numerical ranges of $A$ is the same as in Corollary~\ref{th:hk5}: $\Lambda_1(A)=\conv\{D\cup(-D)\}$, $\Lambda_2(A)=(\conv D)\cap(\conv (-D))$, and $\Lambda_3(A)=\{0\}$.

\begin{figure}[H] 
	\includegraphics[width=8cm]{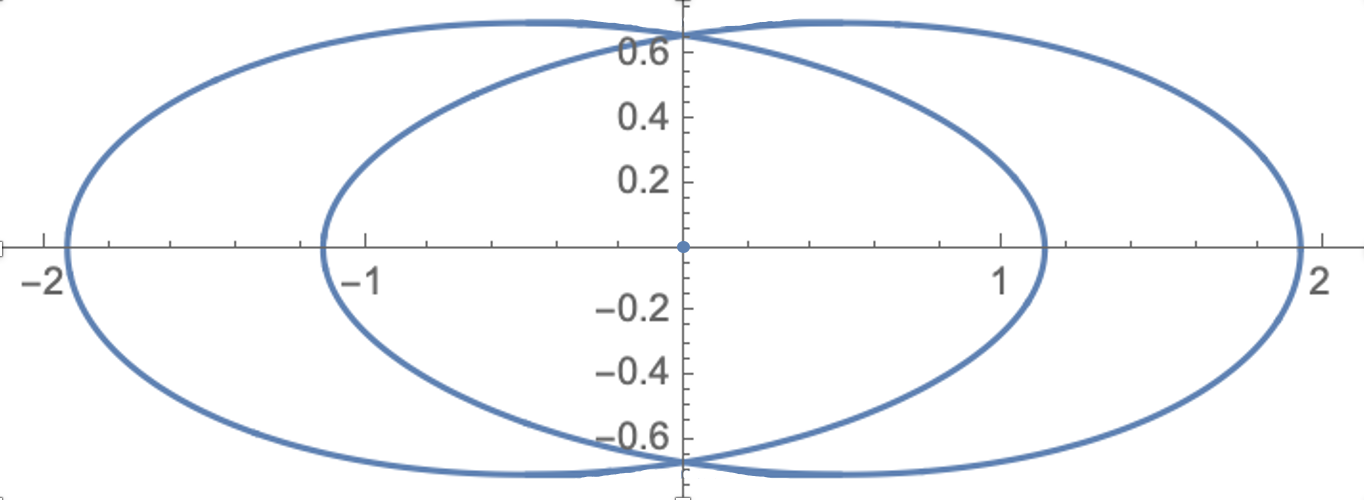}
	\centering
	\caption{Example of $C(A)$ with non-elliptical components, $n=5$} \label{fde5} 
\end{figure}
Observe also that under condition \eqref{xip} $A$ is what in \cite{BPS} was called a 2-periodic reciprocal matrix. Such matrices of dimension $n=1\mod 4$ (in particular, $n=5$) have no elliptical components in $C(A)$ \cite[Theorem 5]{BPS}.
	
To conclude the section, here is an example of $C(A)$ for a matrix $A$ satisfying conditions of Theorem~\ref{th:el5}. It consist of the origin and a pair of displaced ellipses. The rank-2 numerical range is the intersection of the two elliptical  discs.

\begin{figure}[H] 
	\includegraphics[width=8cm]{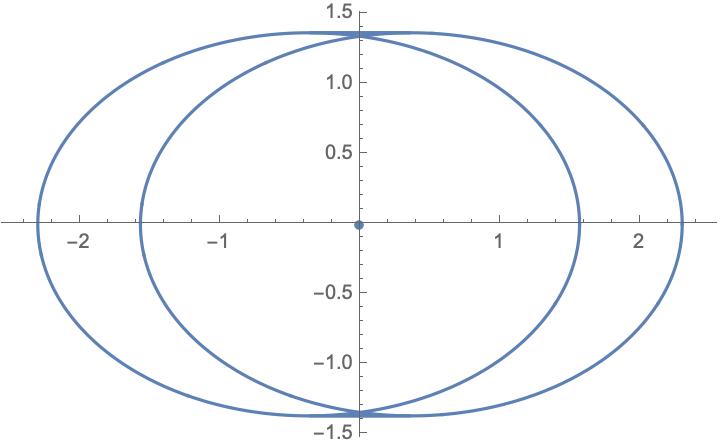}
	\centering
		\caption{Example of $C(A)$ with elliptical components, $n=5$} \label{fde0} 
\end{figure}
\[\{\xi_i\}=\{ 1+\sqrt{3}/{2},0,1,\sqrt{3}/{2}\} \]

\section{Reciprocal 6-by-6 matrices} \label{s:6} 
The Kippenhahn polynomial of $A\in\Rec_6$ is
\begin{multline}
P_6(\zeta,\rho)  = \zeta^3-(\xi_1+\xi_2+\xi_3+\xi_4+\xi_5+\rho)\zeta^2 \\  \label{char6}  +
\left(\textstyle\prod_{j-i>1}\xi_i\xi_j+(3(\xi_1+\xi_5)+2(\xi_2+\xi_3+\xi_4))\rho+6\rho^2\right)\zeta \\ 
-(\xi_1+\rho)(\xi_3+\rho)(\xi_5+\rho).\end{multline} 
In contrast with the $n=4,5$ cases, Proposition~\ref{th:qfac} implies that $C(A)$ may contain both elliptical and non-elliptical components. When this happens, the ellipse contained in $C(A)$ is necessarily centered at the origin, and the respective criterion is delivered by \cite[Theorem 10]{BPSV}. On the other hand, if $C(A)$ contains two ellipses, then it is in fact the union of three ellipses, either all centered at the origin, or one centered at the origin and two symmetric through the ordinate axis. The former case is described by \cite[Theorem 11]{BPSV}, and in our notation (after some additional manipulations) can be restated as follows.
\begin{thm}\label{th:3conel} Let $A\in\Rec_6$. Then for $C(A)$ to consist of three concentric ellipses it is necessary and sufficient that
\eq{3conel}\begin{aligned}
& 2 \xi _1^2+4 \xi _2 \xi _1-2 \xi _3 \xi _1-3 \xi _4 \xi _1-3 \xi _5 \xi _1+\xi _2^2+\xi _4^2+2 \xi _5^2+2 \xi _2 \xi _3\\ & -5 \xi _2 \xi _4
+2 \xi _3 \xi _4-3 \xi _2 \xi _5-2 \xi _3 \xi _5+4 \xi _4 \xi _5 = 0, \\
& 2 \xi _1^2+\xi _2 \xi _1-3 \xi _3 \xi _1+\xi _4 \xi _1-3 \xi _5 \xi _1-\xi _2^2+\xi _3^2-\xi _4^2
+2 \xi _5^2\\ &
+2 \xi _2 \xi _3-2 \xi _2 \xi _4+2 \xi _3 \xi _4+\xi _2 \xi _5-3 \xi _3 \xi _5+\xi _4 \xi _5 = 0, \\ 
& \xi _1^3+2 \xi _2 \xi _1^2+4 \xi _3 \xi _1^2+2 \xi _4 \xi _1^2+3 \xi _5 \xi _1^2-\xi _2^2 \xi _1+3 \xi _3^2 \xi _1-\xi _4^2 \xi _1\\ & 
+3 \xi _5^2 \xi _1+3 \xi _2 \xi _3 \xi _1-2 \xi _2 \xi _4 \xi _1
+3 \xi _3 \xi _4 \xi _1+4 \xi _2 \xi _5 \xi _1-41 \xi _3 \xi _5 \xi _1\\ &
+4 \xi _4 \xi _5 \xi _1-\xi _2^3-\xi _3^3-\xi _4^3+\xi _5^3+2 \xi _2 \xi _3^2-3 \xi _2 \xi _4^2+\xi _3 \xi _4^2\\ &
+2 \xi _2 \xi _5^2+4 \xi _3 \xi _5^2
+2 \xi _4 \xi _5^2+\xi _2^2 \xi _3-3 \xi _2^2 \xi _4+2 \xi _3^2 \xi _4+2 \xi _2 \xi _3 \xi _4\\ & -\xi _2^2 \xi _5+3 \xi _3^2 \xi _5-\xi _4^2 \xi _5+3 \xi _2 \xi _3 \xi _5-2 \xi _2 \xi _4 \xi _5+3 \xi _3 \xi _4 \xi _5=0, \end{aligned} \en
where  $\xi_j$ are defined by \eqref{xij}. 
\end{thm} 

With conditions \eqref{3conel} satisfied, $C(A)=\cup_{j=1}^3 E_j$, where according to \eqref{sA} $E_j$ is an ellipse with the foci $\pm 2\cos\frac{j\pi}{7}$. The lengths of the minor half-axes of $E_j$ are the positive eigenvalues of $\im A$, i.e., positive roots $c_j$ of \eq{im6} P_6(\lambda^2,0)=\lambda^6-\lambda^4\sum_{j=1}^5\xi_j+\lambda^2\textstyle\prod_{j-i>1}\xi_i\xi_j-\xi_1\xi_3\xi_5. \en 
\begin{prop} \label{th:match}In the setting of Theorem~\ref{th:3conel}, the lengths of the axes of $E_j$ agree with the distances between their foci: $c_1\geq c_2\geq c_3$. \end{prop} 
\begin{proof}Let $A(t)$ be a reciprocal matrix obtained from $A$ by replacing its parameters $\xi_i$ with $t\xi_i$ ($t>0$), $i=1,\ldots,5$. Since conditions \eqref{3conel} are homogeneous, all the matrices $A(t)$ will satisfy them along with $A$, and so $C(A(t))=\cup_{j=1}^3 E_j(t)$. Moreover, the polynomial \eqref{im6} is homogeneous in $\lambda^2,\xi_j$, and so the multiset  $\{c_j(t)\}_{j=1}^3$ of its roots is simply $\{\sqrt{t}c_j\}_{j=1}^3$. 

The foci of ellipses $E_j(t)$ do not depend on $t$ while their minor axes depend on $\sqrt{t}$ linearly. Observe that the line segment $I_j:=[-2\cos\frac{j\pi}{7},2\cos\frac{j\pi}{7} ]$ connecting the foci of $E_j(t)$ lies inside it. On the other hand, $E_j(t)$ lies in an arbitrarily small neighborhood of $I_j$ provided that $t$ is chosen small enough. So, for such $t$ the ellipse $E_3(t)$ lies inside $E_2(t)$ which, in turn, lies inside $E_1(t)$. Consequently, $c_1(t)\geq c_2(t)\geq c_3(t)$ for small $t$. But then the inequalities persist for all $t$, since the ellipses $E_j(t)$ remain nested for all $t$  by Corollary~\ref{th:concen}, and so their ordering cannot change. Setting $t=1$ completes the proof. 
\end{proof} 
\begin{cor}Let $A$ be a reciprocal matrix satisfying \eqref{3conel}. Denote by $E_j$ the ellipse with the foci $\pm 2\cos\frac{j\pi}{7}$ and the length of its minor half-axis equal the $j$th (in the non-increasing order) eigenvalue of $\im A$, $j=1,2,3$. Then $\Lambda_j(A)$ is the elliptical disk bounded by $E_j$. \label{th:hrconel} \end{cor} 

We now turn to the case of non-concentric ellipses. As in Section~\ref{s:45}, we first need to figure out when $C(A)$ admits multiple tangent lines. 
\begin{lem} \label{th:hor6}A matrix $A\in\Rec_6$ has the Kippenhahn curve with multiple tangent lines if and only if one of the following (mutually exclusive) conditions holds:\\ {\em (i)} $\xi_1\xi_3\xi_5=0$,\\ {\em (ii)} $\xi_2=\xi_4=0$ while $\xi_1,\xi_3,\xi_5\neq 0$ are not all distinct, \\ {\em (iii)} $\xi_2=0$, $\xi_1\xi_3\xi_4\xi_5\neq 0$,  and \eq{xi14} \xi_1\xi_4=(\xi_1-\xi_5)(\xi_1-\xi_3),\en {\em (iv)} $\xi_4=0$, $\xi_1\xi_2\xi_3\xi_5\neq 0$, and \eq{xi25} \xi_2\xi_5=(\xi_1-\xi_5)(\xi_3-\xi_5).\en  
\end{lem} 
\begin{proof}By Proposition~\ref{th:multan}, for multiple tangent lines to exist it is necessary that at least one $\xi_j$ is equal to zero. So, we just need to consider each of the following situations separately. 
	
{\sl Case 1.} $\xi_1\xi_3\xi_5=0$. Then $\im A$ splits into the direct sum of at least two diagonal blocks each of which is singular. So, $0$ is a multiple eigenvalue, and the abscissa axis is a multiple tangent line. No additional conditions on $\xi_j$ arise in this case.

{\sl Case 2.} $\xi_1\xi_3\xi_5\neq 0,\ \xi_2=\xi_4=0$. The matrix $\im A$ is the direct sum of three 2-by-2 blocks, with the eigenvalues $\pm\sqrt{\xi_j}$ ($j=1,3,5$), respectively. So, $\im A $ has multiple eigenvalues if and only if $\xi_1,\xi_3,\xi_5$ are not all distinct, and (ii) follows. 

{\sl Case 3. $\xi_1\xi_3\xi_5\neq 0$, and exactly one of $\xi_2,\xi_4$ is equal to zero.} Conditions (iii) and (iv) interchange under \eqref{flip}, and so we may without loss of generality suppose that $\xi_2=0,\ \xi_4\neq 0$. 

Similar to \eqref{b} (and with the same notation $\eta_j$) we have $\im A=B_1\oplus B_3$, where 
\[ B_1=\begin{bmatrix}0 & \eta_1 \\ \overline{\eta_1} & 0\end{bmatrix}, \quad B_3=\begin{bmatrix}0 & \eta_3 & 0 & 0 \\ 
\overline{\eta_3} & 0 &  \eta_4 & 0 \\ 0  & \overline{\eta_4} & 0 & \eta_5 \\ 0 & 0 & \overline{\eta_5} & 0 \end{bmatrix}.\]
Since $\eta_3\eta_4\eta_5\neq 0$ along with $\xi_3\xi_4\xi_5$, the eigenvalues of $B_3$ are simple. In order for $\im A$ to have multiple eigenvalues it is therefore necessary and sufficient that the characteristic polynomial $\zeta-\xi_1$ of $B_1$ divides the characteristic polynomial $p(\zeta)=\zeta^2-(\xi_3+\xi_4+\xi_5)\zeta+\xi_3\xi_5$ of $B_3$. In other words, $p(\xi_1)=0$, which is exactly \eqref{xi14}. This takes care of (iii) and, after invoking \eqref{flip}, also (iv). 
\end{proof} 
\begin{thm}\label{th:el6}Let $A\in\Rec_6$. Then $C(A)$ consists of three ellipses exactly one of which is centered at the origin if and only if one of the following conditions holds:
\eq{de1} \xi_2=\xi_4=2\xi_1\cos\frac{2\pi}{7},\ \xi_3=0,\ \xi_5=\xi_1,\en 
\eq{de2} \xi_3=\xi_5=k\xi_1,\ \xi_2=0,\ \xi_4=(k-1)^2\xi_1,  \en 
or 
\eq{de3} \xi_1=\xi_3=k\xi_5,\ \xi_2=(k-1)^2\xi_5,\ \xi_4=0, \en 
with the coefficient $k$ in \eqref{de2},\eqref{de3} taking the values $2\cos\frac{\pi}{7}$, $2\cos\frac{3\pi}{7}$.
\end{thm}
\begin{proof}By Proposition~\ref{th:facrec}, \eq{3el} C(A)=E_0\cup E\cup (-E), \en where $E_0$ ($E$) is an ellipse centered at the origin (resp., some $p>0$) if and only if the polynomial \eqref{char6} factors as 
\eq{fac6}
\left(\zeta -\left(c_0^2+\rho  X_0^2\right)\right) \left(\zeta ^2-2 \zeta  \left(c^2+\rho  \left(p^2+X^2\right)\right)+\left(c^2+\rho  \left(X^2-p^2\right)\right)^2\right).   \en 
Recall also that $X_0$ and $c_0$ ($X$ and $c$) are half the distance between the foci and half the length of the minor axis of $E_0$ (resp. $E$).  	

Equating the respective coefficients of \eqref{char6} and \eqref{fac6}, we arrive at the system of nine equations which for our purposes is convenient to group into three subsystems \eqref{Xp}--\eqref{xiXp}, with three equations in each:   	
\eq{Xp}  \begin{cases} X_0^2(X^2-p^2)^2=1 \hspace{6cm} \\
	(X^2-p^2)^2+2X_0^2(X^2+p^2)=6\\
	X_0^2+2(X^2+p^2)=5, 
\end{cases}\en 
\eq{xic1}\begin{cases}		\xi _1 \xi _3 \xi _5=c^4 c_0^2 \hfill  \\
		\xi_1+\xi _2+\xi _3+\xi_4+\xi _5=2 c^2+c_0^2\\
		\xi _1 \xi _3+\xi _1 \xi _4+\xi _2 \xi _4+\xi _1 \xi _5+\xi _2 \xi _5+\xi _3 \xi _5=2 c^2 c_0^2+c^4, \end{cases} \en  
\eq{xiXp}\begin{cases}
		 \xi _1 \xi _3+\xi _1 \xi _5+\xi _3 \xi _5=c^4 X_0^2+2 c^2c_0^2 (X^2-p^2)\\
		\xi_1+\xi _3+\xi _5=2 c^2 X_0^2(X^2-p^2)+c_0^2(X^2-p^2)^2\\
		3(\xi_1+\xi_5)+2(\xi _2+\xi _3+\xi_4)= 2 c^2 (X_0^2+X^2-p^2)+2c_0^2(X^2+p^2). \end{cases} \en 
The first of them, \eqref{Xp}, contains only the variables $X_0,X,p$ and (taking into consideration that they are positive, due to their geometrical meaning) has exactly six solutions, as per the table below: 
\begin{center}
	\begin{tabular}{ | c | c | c | c |} 
		\hline
		& $X_0$ & $X$  & $p$ \\
	 \hline
	(i)	& $2\cos\frac{\pi}{7}$ & $\cos\frac{2\pi}{7}-\cos\frac{3\pi}{7}$ & $\cos\frac{2\pi}{7}+\cos\frac{3\pi}{7}$ \\ 
		\hline
	(ii) & 	$2\cos\frac{\pi}{7}$ & $\cos\frac{2\pi}{7}+\cos\frac{3\pi}{7}$ & $\cos\frac{2\pi}{7}-\cos\frac{3\pi}{7}$ \\ 
		\hline
	(iii) &	$2\cos\frac{2\pi}{7}$ & $\cos\frac{\pi}{7}-\cos\frac{3\pi}{7}$ & $\cos\frac{\pi}{7}+\cos\frac{3\pi}{7}$ \\ 
		\hline
(iv) & $2\cos\frac{2\pi}{7}$ & $\cos\frac{\pi}{7}+\cos\frac{3\pi}{7}$ &  $\cos\frac{\pi}{7}-\cos\frac{3\pi}{7}$ \\
\hline 
(v) & $2\cos\frac{3\pi}{7}$ & $\cos\frac{\pi}{7}-\cos\frac{2\pi}{7}$ & $\cos\frac{\pi}{7}+\cos\frac{2\pi}{7}$ \\
\hline 
(vi) & $2\cos\frac{3\pi}{7}$ & $\cos\frac{\pi}{7}+\cos\frac{2\pi}{7}$ &  $\cos\frac{\pi}{7}-\cos\frac{2\pi}{7}$ \\
\hline 	
\end{tabular}
\end{center} 
To move forward, recall that Lemma~\ref{th:hor6} delivers necessary conditions for $C(A)$ to be as in \eqref{3el}. We will therefore consider separately Cases~1--3, as introduced in the proof of that Lemma. 

{\sl Case 1.} $\xi_1\xi_3\xi_5=0$. Let first $\xi_3=0$. Solving \eqref{xic1}, we conclude that then either (a)  $c=0$, in which case also $\xi_1=\xi_2=0$ or $\xi_4=\xi_5=0$, or (b) $c_0=0$, $\xi_1+\xi_2=\xi_4+\xi_5=c^2$. Subcase (a) is incompatible with all the solutions of \eqref{xiXp}. In turn, under (b) system \eqref{xiXp} is consistent if and only if $X_0,X,p$ are given by (vi), and the solution is the one-parameter family \eqref{de1}.

Suppose now that $\xi_3\neq 0$ but $\xi_1\xi_5=0$. Due to the symmetry \eqref{flip}, it suffices to consider $\xi_5=0$. From \eqref{xic1} we then conclude that either 
\eq{xi51} c=0,\ \xi_1+\xi_2+\xi_3+\xi_4=c_0^2,\ \xi_1=\xi_2\xi_4=0, \en or \eq{xi52} c_0=0,\ \xi_1+\xi_2+\xi_3+\xi_4=2c^2,\ \xi_1\xi_3+\xi_1\xi_4+\xi_2\xi_4=c^4.\en  
Plugging \eqref{xi51} into the last equation of \eqref{xiXp} yields $c_0^2=c_0^2(X^2+p^2)$. Since for al the solutions (i)--(vi) of \eqref{Xp}, $X^2+p^2\neq 1$, this is an impossibility. 

A more involved (but still direct) computation shows that \eqref{xi52} also is in conflict with the system \eqref{xiXp}. So, the subcase $\xi_5=0$ (and thus also $\xi_1=0$) does not materialize. 

{\sl Case 2.} $\xi_1\xi_3\xi_5\neq 0,\ \xi_2=\xi_4=0$. Under these conditions, \eqref{xic1} means that $\{\xi_1,\xi_3,\xi_5\}$ and $\{c^2,c^2,c_0^2\}$ coincide as multisets. The first two equations of \eqref{xiXp} can then be rewritten as the homogeneous system (with respect to the unknowns $c^2,c_0^2$) with the matrix of coefficients \[ \begin{bmatrix} X_0^2-1 & 2(X^2-p^2-1) \\
	2(X_0^2(X^2-p^2)-1) & (X^2-p^2)^2-1 \end{bmatrix}. \] 
A straightforward computation shows that this matrix is invertible for all the choices (i)--(vi) of solutions to \eqref{Xp}. This is a contradiction with $\xi_j$ being different from zero for $j=1,3,5$. So, no solutions to \eqref{Xp}--\eqref{xiXp} emerge in this case. 

{\sl Case 3. $\xi_1\xi_3\xi_5\neq 0$, and exactly one of $\xi_2,\xi_4$ is equal to zero.} 

Let $\xi_2=0$, $\xi_4\neq 0$. As was shown in the respective part of the proof of Lemma~\ref{th:hor6}, $\xi_1$ is then the multiple eigenvalue of $\im A$ if and only if \eqref{xi14} holds. Under this condition we therefore have \eq{xi42} c_0^2=\xi_1,\quad c^2=\xi_3+\xi_4+\xi_5-\xi_1. \en 
Given $\xi_2=0$ and \eqref{xi42}, \eqref{xic1} is actually equivalent to \eqref{xi14}. Yet another computation shows that \eqref{de2} is the solution of \eqref{xiXp} satisfying also \eqref{xi14} and \eqref{Xp}. More specifically, the solution with $k=2\cos\frac{3\pi}{7}$ corresponds to $X_0,X,p$ given by (ii), $k=2\cos\frac{\pi}{7}$ corresponds to (vi), and (i), (iii)--(v) yield no solutions. 

Invoking \eqref{flip}, we immediately conclude from here that \eqref{de3} is the solution of \eqref{Xp}--\eqref{xiXp} corresponding to the case $\xi_2\neq 0$, $\xi_4=0$. We have exhausted all the possibilites, so the proof is now complete. 
\end{proof} 

Retracing the proof of Theorem~\ref{th:el6}, a more detailed description of the configuration \eqref{3el} can be obtained. Namely:

The length of the minor half-axis of $\pm E$ is $\sqrt{\xi_1+\xi_2}$ ($=\sqrt{\xi_4+\xi_5}$) in case \eqref{de1}, $\sqrt{\xi_3+\xi_4+\xi_5-\xi_1}=k\sqrt{\xi_1}$ in case \eqref{de2}, and $\sqrt{\xi_1+\xi_2+\xi_3-\xi_5}=k\sqrt{\xi_5}$ in case \eqref{de3}.
For the ``central'' ellipse $E_0$ the respective values are 0 (so that in case \eqref{de1} $E_0$ degenerates into the doubleton of its foci), $\sqrt{\xi_1}$ and $\sqrt{\xi_5}$. 

In the setting of \eqref{de1}, or \eqref{de2}, \eqref{de3} with $k=2\cos{\frac{3\pi}{7}}$, the foci of the ``right'' ellipse $E$ are $2\cos\frac{\pi}{7}$ and $-2\cos\frac{2\pi}{7}$, while $E_0$ has the foci $\pm 2\cos\frac{3\pi}{7}$. On the other hand, for \eqref{de2}, \eqref{de3} with $k=2\cos{\frac{\pi}{7}}$, the foci of $E$ are $2\cos\frac{2\pi}{7}$ and $-2\cos\frac{3\pi}{7}$, while $E_0$ has the foci $\pm 2\cos\frac{\pi}{7}$.

Observe that $E_0$ either lies inside the intersection of $E$ with $-E$ or contains their union. A moment's thought reveals that no other configurations are possible because of Proposition~\ref{th:multan}. Also, not only the ellipses $\pm E$ overlap (in agreement with part (i) of Corollary~\ref{th:Aj1}), but moreover, the sets of their foci interlace. 

Here are some numerical examples corresponding to each of the cases described above.

\begin{figure}[H]
	\includegraphics[width=8cm]{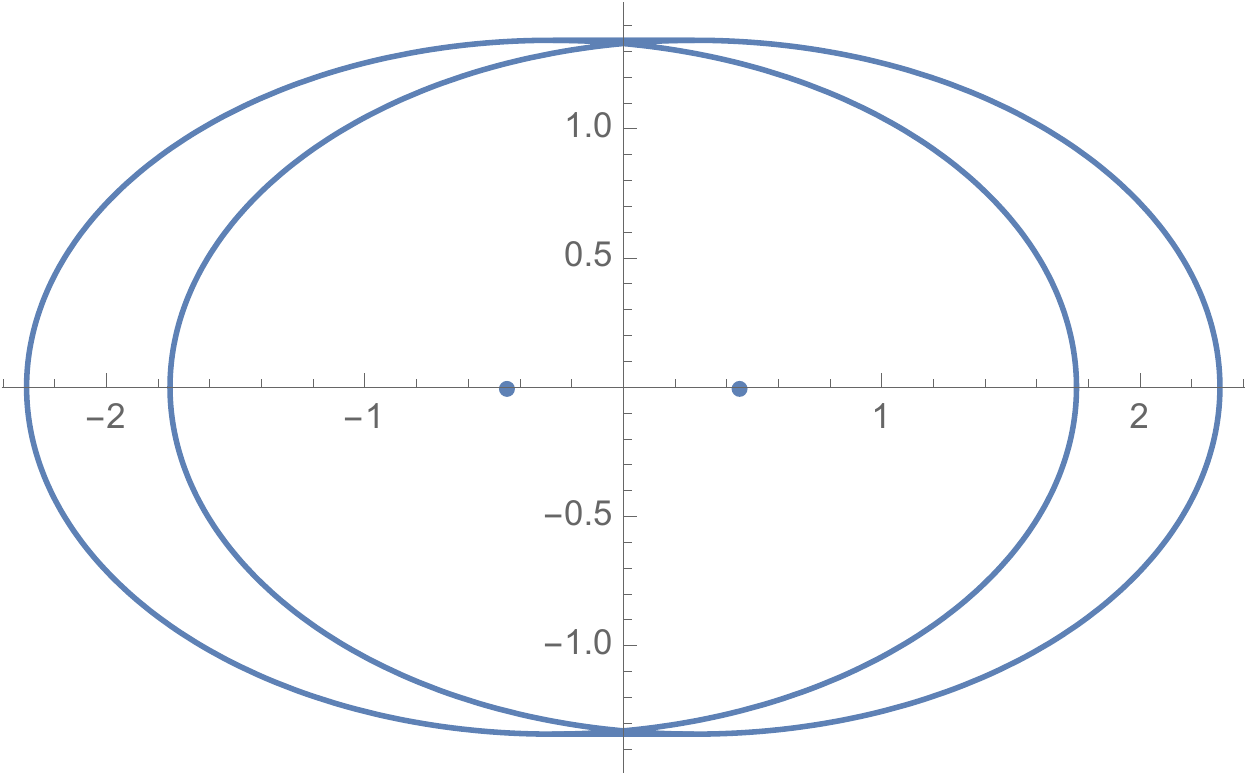}
	\centering
	\caption{Example of C(A) with degenerated central ellipse} \label{fde1} 
\end{figure}
\[\{\xi_i\} = \{0.801938,1.,0.,1.,0.801938\}\]

\begin{figure}[H]
	\includegraphics[width=8cm]{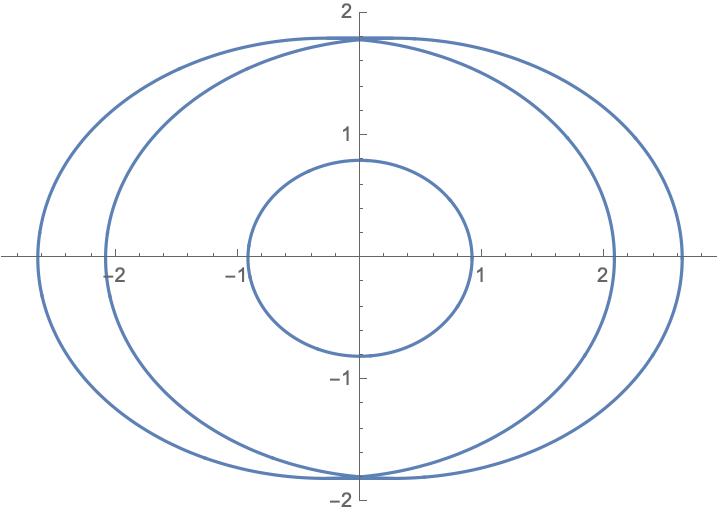}
	\centering
	\caption{Example of C(A) with inner central ellipse}
\end{figure}
\[\{\xi_i\} = \{1.44504,1.,1.44504,0.,3.24698\}\]

\begin{figure}[H]
	\includegraphics[width=8cm]{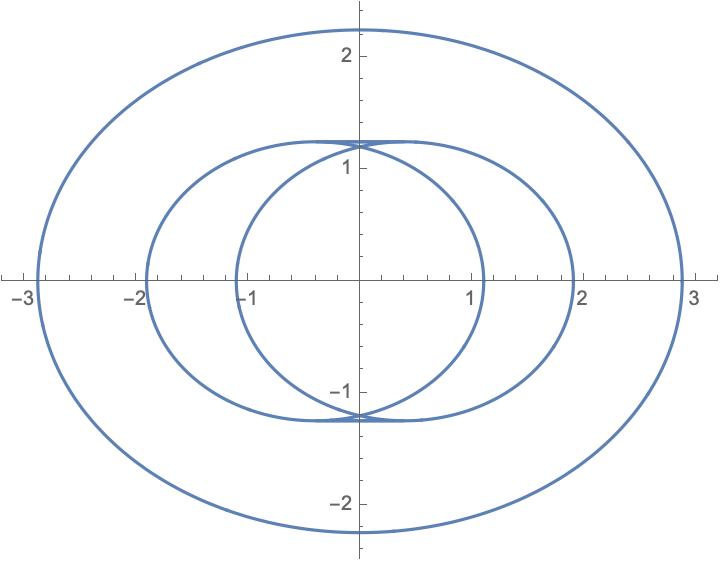}
	\centering
	\caption{Example of C(A) with outer central ellipse}
\end{figure}
\[\{\xi_i\} = \{2.80194,1.,2.80194,0.,1.55496\}\]

Based on the above description of $C(A)$, we arrive at the following result concerning the rank-$j$ numerical ranges of $A$ for $j=1,2,3$ (note that $\Lambda_j(A)=\emptyset$ for $j\geq 4$ due to Corollary~\ref{th:bigk}).
\begin{cor}\label{th:hr3el}In the setting of of Theorem~\ref{th:el6}, the rank-$j$ numerical ranges of $A$ are as described by the following table:
	\begin{center}
		\begin{tabular}{ | c | c | c | } 
			\hline
			& \eqref{de1} or \eqref{de2}, \eqref{de3} & \eqref{de2}, \eqref{de3} \\
			&  with $k = 2\cos{\frac{3\pi}{7}}$ & with $k = 2\cos{\frac{\pi}{7}}$\\ \hline
			$\conv E_0$ & $\Lambda_3$ & $\Lambda_1$ \\ 
			\hline
			$\conv(E \cup -E)$ & $\Lambda_1$ & $\Lambda_2$ \\ 
			\hline
			$\conv E \cap \conv(-E)$ & $\Lambda_2$ & $\Lambda_3$ \\ 
			\hline
		\end{tabular}
\end{center}\end{cor} 
Note that under condition \eqref{de1} $\conv E_0$ is the line segment $I_3$ while under \eqref{de2} or \eqref{de3} it is a  non-degenerate elliptical disk. 
\newpage 
\providecommand{\bysame}{\leavevmode\hbox to3em{\hrulefill}\thinspace}
\providecommand{\MR}{\relax\ifhmode\unskip\space\fi MR }
\providecommand{\MRhref}[2]{%
	\href{http://www.ams.org/mathscinet-getitem?mr=#1}{#2}
}
\providecommand{\href}[2]{#2}

\end{document}